\documentclass[12pt,reqno,oneside]{amsart}
\usepackage{amsmath}
\usepackage{amssymb}
\usepackage{amsthm}
\usepackage{graphics}
\usepackage{eucal}
\usepackage{mathrsfs}
\usepackage[pdftex,bookmarks,colorlinks,breaklinks]{hyperref}  
\newtheorem{theorem}{Theorem}[section]

\newtheorem{proposition}{Proposition}[section]

\numberwithin{equation}{section}

\newcommand{\be}{\mathbf{e}}
\newcommand{\eo}{\mathbf{e}_{0}}

\newcommand{\ea}{\mathbf{e}_{1}}
\newcommand{\ei}{\mathbf{e}_{i}}

\newcommand{\en}{\mathbf{e}_{n}}

\newcommand{\eej}{\mathbf{e}_{J}}

\newcommand{\esj}{\mathbf{e}_{*j}}

\newcommand{\esa}{\mathbf{e}_{*1}}

\newcommand{\esd}{\mathbf{e}_{*d}}

\newcommand{\diag}{\operatorname{diag}}

\newcommand{\Span}{\operatorname{span}}
\newcommand{\GL}{\operatorname{GL}}
\newcommand{\Mat}{\operatorname{Mat}}
\newcommand{\bbR}{\mathbb{R}}

\newcommand{\bbZ}{\mathbb{Z}}

\newcommand{\bq}{\mathbf{q}}
\newcommand{\bx}{\mathbf{x}}

\newcommand{\scH}{\mathcal{H}}

\newcommand{\bp}{\mathbf{p}}

\newcommand{\La}{\Lambda}

\newcommand{\Ga}{\Gamma}
\newcommand{\bv}{\mathbf{v}}
\newcommand{\bw}{\mathbf{w}}
\newcommand{\f}{\mathbf{f}}

\newcommand{\pfij}{\frac{\partial f_j (x)}{\partial x_i}}

\newcommand{\Id}{\operatorname{Id}}

\newcommand{\bc}{\mathbf{c}}
\newcommand{\pb}{\pi_{\bullet}}
\newcommand{\oj}{\omega_{j}}

\newcommand{\al}{\alpha}

\newcommand{\cW}{\mathcal{W}}

\begin{document}
\title[Quantitative Diophantine approximation]{Quantitative Diophantine approximation on affine subspaces}
\begin{abstract}
Recently, Adiceam et.al. \cite{ABLVZ} proved a quantitative version of the convergence case of the Khintchine-Groshev theorem for nondegenerate manifolds, motivated by applications to interference alignment. In the present paper, we obtain analogues of their results for affine subspaces. 
\end{abstract}


\author{Arijit Ganguly and Anish Ghosh}
\thanks {Ghosh is supported by an ISF-UGC grant}
\address{School of Mathematics,
Tata Institute of Fundamental Research, Mumbai, India 400005}
\email{arimath@math.tifr.res.in, ghosh@math.tifr.res.in}

\maketitle

\section{Introduction}

The theory of Diophantine approximation on manifolds has seen significant  advances in recent years. This subject is mainly concerned with the question: under which conditions do proper subsets of $\mathbb{R}^n$ inherit Diophantine properties which are generic for $\mathbb{R}^n$ with respect to Lebesgue measure? A simple example of such a generic Diophantine property is provided by the classical  Khintchine-Groshev theorem. Let $\psi : \bbR_{+} \cup \{0\} \to \bbR_{+} \cup \{0\}$ be a non-increasing function and consider the set of $\psi$-approximable vectors, namely $x \in \mathbb{R}^n$ for which there exist infinitely many $\bq \in \bbZ^{n}$ such that
\begin{equation}\label{preKG}
|p +  \bx \cdot \bq| < \psi(\|\bq\|^n)
\end{equation}
\noindent for some $p \in \bbZ$. We will use $|~|$ for both the Lebesgue measure of a measurable subset of $\bbR^n$ as well as the absolute value of a real number.  The Khintchine-Groshev Theorem (\cite{Khintchine}, \cite{Groshev}, \cite{Dodson}) states that the set of $\psi$-approximable vectors is a null (reap. co-null) set in terms of Lebesgue measure, according as the sum

\begin{equation} 
 \sum_{k = 1}^{\infty} \psi(k) 
\end{equation}

\noindent converges or diverges. Let $U$ be an open subset of $\bbR^d$ and let $\f : U \to \bbR^n$ be a differentiable map. Then $\f$ is said to be \emph{nondegenerate} at $x \in U$ if $\bbR^n$ is spanned by the
partial derivatives of $\f$ at $x$ of order up to $l$ for some $l$, and nondegenerate if it is nondegenerate at almost every point of $U$. Nondegenerate manifolds, i.e. manifolds parametrised by nondegenerate maps, inherit many generic Diophantine properties from ambient Euclidean space. For instance, in an influential paper  D. Kleinbock and G. Margulis \cite{KM} resolved a long standing conjecture of Sprind\v{z}uk by showing that nondegenerate maps are extremal, i.e. almost every point on such a manifold is \emph{not very well approximable}. 

Subsequently, V. Bernik, D. Kleinbock and G. Margulis \cite{BKM} established the convergence case of Khintchine's theorem for nondegenerate manifolds. This result was independently established by V. Beresnevich \cite{Ber}. In fact, both \cite{KM} and \cite{BKM} prove multiplicative versions of these results. In \cite{ABLVZ}, Adiceam et.al. have recently proved an interesting quantitative improvement of the convergence Khintchine theorem for nondegenerate manifolds. Their motivation comes from electronics, more precisely the study of interference alignment.

At the opposite end of the spectrum from nondegenerate manifolds lie affine subspaces. The study of Diophantine approximation on affine subspaces and their submanifolds was systematically initiated in the works \cite{Kleinbock-extremal, Kleinbock-exponent} of D. Kleinbock and has seen recent progress. Since arbitrary affine subspaces cannot be expected to inherit generic Diophantine properties, the interesting question of finding necessary and sufficient conditions on affine subspaces to ensure inheritance of a given property plays a key role in investigations. We refer the reader to the recent survey \cite{G-handbook} for a comprehensive discussion as well as references. 
 
In this paper, we undertake the study of the refined, quantitative, version of the Khintchine-Groshev theorem from \cite{ABLVZ} in the context of affine subspaces. We provide a sufficient condition for an affine subspace to satisfy such a theorem. This condition is introduced in the next subsection after which we state the main result of the paper. In addition to the interest in this problem from the Diophantine point of view, it is possible that the result proved here could have applications in interference alignment. This is explained in \cite{ABLVZ}, and indeed the example presented there, concerns a line!

\subsection{Diophantine exponents of matrices}\label{sec:exp}
Let $\scH$ be an $s$ dimensional affine subspace of $\bbR^n$. We can permute variables and assume that $\scH$ is of the 
form $\{ (\mathbf{x},\mathbf{x} A' + \mathbf{a}_0): \mathbf{x}\in \mathbb{R}^s\}$ where $\mathbf{a}_0 \in \mathbb{R}^{n-s}$ and $A' \in \Mat_{s \times n-s}(\bbR)$. Denoting 
the matrix  $\begin{pmatrix}\mathbf{a}_0\\A' \end{pmatrix}$ by $A$, we can rewrite the parametrization as 
\begin{equation}\label{defaffine}
\bx \mapsto (\bx,  \tilde{\bx}A)\,\text{where }\tilde{\bx} = (1, \bx).
\end{equation}

\noindent The Diophantine exponent $\omega(A)$ of a matrix $A \in \Mat_{m \times n}(\bbR)$ is defined to be the 
supremum of $v > 0$ for which there are infinitely many $\bq \in \bbZ^n$ such that
\begin{equation}\label{defexponent}
\|A \bq + \bp\| < \|\bq\|^{-v}
\end{equation}
\noindent for some $\bp \in \bbZ^m$. It is well known that $n/m \leq \omega(A) \leq \infty$ 
for all $A \in \Mat_{m \times n}(\bbR)$ and that $\omega(A) = n/m$ for Lebesgue almost 
every $A$. We now introduce the higher Diophantine 
exponents of $A$ as defined by Kleinbock in \cite{Kleinbock-exponent}. For $A \in \Mat_{s+1 \times n-s}(\bbR)$, we set 
\begin{equation}\label{defR}
R_A = \begin{pmatrix}\Id_{s+1} & A \end{pmatrix}.
\end{equation}Let $\eo, \dots, \en$ denote the standard basis of $\bbR^{n+1}$ and set 
\begin{equation}
W_{i \to j} = \Span\{\be_i, \dots, \be_j\}.
\end{equation}
\noindent Let $\bw \in \bigwedge^{j}(W_{0 \to n})$ represent a discrete subgroup $\Gamma$ of $\bbZ^{n+1}$. Define 
the map $\bc : \bigwedge^{j}(W_{0 \to n}) \to (\bigwedge^{j-1}(W_{1 \to n}))^{n+1}$ by
\begin{equation}\label{defc}
\bc(\bw)_i = \sum_{\substack{J \subset \{1, \dots, n\}\\ \#J = j-1}} \langle \be_i \wedge \eej, \bw\rangle \eej
\end{equation} 
\noindent and let $\pb$ denote the projection $\bigwedge(W_{0 \to n}) \to \bigwedge(W_{s+1 \to n})$. For each $j = 1, \dots, n-s$, define
\begin{equation}\label{defexponenthigher}
\oj(A) = \sup\left\{ v\left| \aligned \exists\, 
\bw \in  \bigwedge^{j}(\bbZ^{n+1}) \text{ with arbitrary large
   } \|\pb(\bw)\|  \\
\text{   such that   }\|R_{A}
\bc(\bw)\|
< \|\pb(\bw)\|^{-\frac {v+1-{j}}{j}}\ \ 
\endaligned\right.\right\}.
\end{equation}
\noindent It is shown in Lemma 5.3 of \cite{Kleinbock-exponent} that $\omega_1(A) = \omega(A)$ thereby justifying 
the terminology.

\subsection{Main Theorem}
Let $\psi$ be an approximation function with $\psi(x)\,\leq\,\frac{1}{x}$ for all $x\in \mathbb{R}$. Assume that 
\begin{equation}\label{eqn:conv1}
 \displaystyle \sum_{k=1} ^{\infty}  \psi(k)\,\textless \,\infty\,.
\end{equation}
Since $\psi$ is assumed  non-increasing, it is easy to see that condition (\ref{eqn:conv1}) is equivalent to saying that 
\begin{equation}\label{eqn:conv2}
 \textstyle\sum_{\psi} := \displaystyle\sum_{\bq\in \mathbb{Z}^n\setminus \{\mathbf{0}\}} ^{\infty}  \psi(||\bq||^n)\,\textless \,\infty\,.
\end{equation}

We consider $\psi$-approximable points on affine subspaces, namely solutions to the inequality
\begin{equation}\label{reform}
| (\bx, \tilde{\bx}A) \cdot \bq + p| < \psi(\|\bq\|^n).
\end{equation}
As a corollary of (Theorem 1.2, \cite{G-mani}),  we see that for any open ball $U$ in $\mathbb{R}^s$, the measure of the set
\begin{equation}\label{eqn:BC}
\{\mathbf{x}\in U: \exists~p\in \mathbb{Z} \text{ such that (\ref{reform}) holds for \text{infinitely many} }\bq\in \mathbb{Z}^n\} 
 \end{equation}
is zero,  provided  (\ref{eqn:conv2}) holds and
 \begin{equation}\label{diocond}
\oj(A) < n~\text{for every}~j = 1, \dots, n-s. 
\end{equation} 
Thus for almost all $\bx\in U$, there exists a constant $\kappa \,\textgreater\,0$ such that 
\begin{equation}\label{eqn:kappa1}
 |(\bx,  \tilde{\bx}A)\bq + p|\,\geq\, \kappa \,\psi(\|\bq\|^n)\text{ for all }p\in \mathbb{Z}, \bq \in \mathbb{Z}^n\setminus \{\mathbf{0}\}\,.
\end{equation}

\noindent 
Consider the set
\begin{equation}\label{eqn:defB}
 \mathcal{B}(U, \psi, \kappa):=\{\mathbf{x}\in U:\text {(\ref{eqn:kappa1}) holds} \}\,.  \end{equation}
Our aim is to investigate the dependence between $\kappa$ and the size of the set (\ref{eqn:defB}). Our main Theorem is
\begin{theorem}\label{main}
Let $\scH$  be an $s$-dimensional affine subspace parametrized as in (\ref{defaffine}). Assume that
\begin{equation}\label{diocond}
\oj(A) < n~\text{for every}~j = 1, \dots, n-s. 
\end{equation} Consider a non-increasing approximation function $\psi$ such 
that $\psi(x)\,\leq\,\frac{1}{x}$ for all $x\in \mathbb{R}$ and assume that that (\ref{eqn:conv1}) holds. 
Fix an open ball $U\in \mathbb{R}^s$. Then there exist two explicitly computable constants $K_0$ and $K_1$, depending on $s,n, U\text{ and }A$ only,
with the following property: \\

\indent for any $\xi \in (0,1)$,  
\begin{equation}\label{eqn:goal}
|\mathcal{B}(U, \psi, \kappa)|\geq (1-\xi)|U|,
\end{equation} holds with 
\begin{equation}\label{eqn:kappa2}
 \kappa< \displaystyle \min \left \{1,\frac{\xi}{2 K_s\textstyle\sum_{\psi}},\frac{r}{2^{n-\frac{3}{2}}\sqrt{ns}}, 
 \left(\frac{\xi}{2K_0K_1}\right)^{s(n+1)}
 \right\}.
\end{equation} Where \[K_{s}:=\frac{4^{2s+1}s^{s/2}N_s}{V_s}\,,\] $V_s$ is the volume of the $s$-dimensional euclidean unit
ball and $N_s$ denotes the Besicovitch covering constant of $\mathbb{R}^s$.  

\end{theorem}

Remarks:
\begin{enumerate}
\item Although we have not pursued it here, it is plausible that Theorem \ref{main} is also true for nondegenerate subaminfolds of affine subspaces under the same condition, i.e. (\ref{diocond}).
\item We will follow the general strategy of of \cite{ABLVZ} to prove Theorem \ref{main}, indeed this can be traced back to the work of Bernik, Kleinbock and Margulis \cite{BKM}.  The proof splits into two separate cases, the `big gradient' and `small gradient'. Most of this paper is devoted to the latter case, and involves nondivergence estimates for polynomial like flows on the space of unimodular lattices. 
\item It is worthwhile considering the case where $\scH$ is a hyperplane namely an $n-1$ dimensional subspace of $\mathbb{R}^n$. In this case, $A$ is an $n \times 1$ matrix and the condition (\ref{diocond}) takes a particularly simple form, namely that for some  $\delta > 0$,
\begin{equation}\nonumber
\max_{i}|p_i + a_iq| > |q|^{-n+\delta}
\end{equation}
for every $p \in \mathbb{Z}^n$, and all but finitely many $q \in \mathbb{Z}$. 
\end{enumerate}

\section*{Acknowledgements} Part of this work was done when the second named author was visiting Peking University. He thanks J. An for his hospitality.

\section{The Gradient Division}\label{gradientdivision}
For $\kappa \,\textgreater\, 0\text{ and }\mathbf {q} \in \mathbb{Z}^n\setminus \{\mathbf{0}\}$, we define 
\[\mathcal {L}(\mathbf {q}):= \{ \mathbf{x}\in U: |p + (\bx,  \tilde{\bx}A)\bq|\,\textless\, \kappa \,\psi(\|\bq\|^n)\text{ for some }p\in \mathbb{Z}\}\,.\]
As $\displaystyle \bigcup_{\mathbf{q}\in \mathbb{Z}^n \setminus \{\mathbf{0}\}} \mathcal {L}(\mathbf {q})= U\setminus \mathcal{B}(U, \psi, \kappa)$, 
it suffices to prove that 
$$\displaystyle \big |\bigcup_{\mathbf{q}\in \mathbb{Z}^n \setminus \{\mathbf{0}\}} \mathcal {L}(\mathbf {q})\big |\,\textless\,\xi|U|\,.$$ 

It is traditional to approach Khintchine-Groshev type theorems by
separately considering the case when we have a `large derivative'
and the case when we do not. 
We are thus interested in the
cases where
$\nabla (\bx,  \tilde{\bx}A)\cdot \bq =[Id_s\,\,A']\bq$, where $A'$ is as introduced in the beginning of \S\ref{sec:exp},  gets big or small. Let
\begin{equation}
\mathcal {L}_{small}(\mathbf {q}) = \left\{x \in \mathcal {L}(\mathbf {q})~:~ 
\|\nabla (\bx,  \tilde{\bx}A)\cdot \bq\|<\frac{\sqrt{ns \|\bq\|}}{2r} \right \}\end{equation}
\noindent where $r$ is the radius of $U$ and $\mathcal {L}_{large}(\mathbf {q}) = \mathcal {L}(\mathbf {q}) 
\backslash \mathcal {L}_{small}(\mathbf {q})$. We will prove that for $\kappa$ 
given by (\ref{eqn:kappa2}),
\begin{equation}\label{Borcant1}
\sum_{\bq \in \bbZ^{n}} |\mathcal {L}_{large}(\mathbf {q})| \,\leq \,\frac{\xi}{2}|U|
\end{equation}
\noindent and
\begin{equation}\label{Borcant2} 
\big |\bigcup_{\bq \in \bbZ^{n}}\mathcal {L}_{small}(\mathbf {q})\big |\, <\, \frac{\xi}{2}|U|.
\end{equation}
\section{Estimating the measure of $\mathcal {L}_{large}(\mathbf {q})$}
In this section, we will establish (\ref{Borcant1}). The proof of this follows immediately from 
\begin{proposition}\cite[Theorem 4]{ABLVZ}\label{prop:large}
Let $U\subseteq \mathbb{R}^s$ be a ball of radius $r$ and  $\f \in C^2 (2U)$ where $2U$ is the ball with the same center as U and radius $2r$.
Set \begin{equation}\label{defL*}
L^*:=\sup_{|\beta|=2,\,\mathbf{x}\in 2U}\|\partial_{\beta}\f(\bx)\|
\end{equation}and 
\begin{equation}\label{defL}
 L:=\max  \left \{L^*, \frac{1}{4r^2}\right \}.
\end{equation} Then for every  $\delta'\,\textgreater\,0$ and every $\bq\in \mathbb{Z}^n\setminus \{\mathbf{0}\}$,  the set
of all $\mathbf{x}\in U$ such that $|p+\nabla \f(\bx)\bq|\,\textless\,\delta'$ for some $p\in \mathbb{Z}$ and 
\begin{equation}
\|\nabla \f(\bx)\bq\|\,\geq\,\sqrt{nsL\|\bq\|}
\end{equation}
\noindent has measure at most $K_{s}\delta'|U|$.
\end{proposition}

\noindent The proof  of Proposition \ref{prop:large} is done by applying \cite[Lemma $2.2$]{BKM} appropriately.\\

To prove (\ref{Borcant1}) from Proposition \ref{prop:large}, we take $\f(\bx)=(\mathbf{x},\tilde{\mathbf{x}}A)$ 
and $\delta'=\kappa\, \psi(\|\bq\|^n)$. Clearly $L^*=0$ and $L=\frac{1}{4r^2}$. Hence by Proposition \ref{prop:large}, we get that 
\[|\mathcal {L}_{large}(\mathbf {q})|\,\leq\,K_s \kappa \,\psi(\|\bq\|^n)|U|\,,\] and thus, taking $\kappa \leq\frac{\xi}{2 K_s\sum_{\psi}}$,  
\[\sum_{\bq \in \bbZ^{n}} |\mathcal {L}_{large}(\mathbf {q})| \,\leq \,K_s \kappa 
\textstyle\sum_{\psi} |U|\,\leq\, \displaystyle \frac{\xi}{2}|U|\,.\,\,\, \Box\]
To estimate $|\mathcal {L}_{small}(\mathbf {q})|$, we shall employ dynamical tools. To begin with, we need to recall a few elementary
properties of `good functions' which will be discussed in the following
section.

\section{$(C, \alpha)$-good functions}
\noindent Let $C$ and $\al$ be positive numbers and $V$ be a subset of
$\bbR^s$. A function $f : V \to \bbR$ is said to be $(C,\al)$-good on
$V$ if for any open ball $B \subseteq V$,~and
for any $\varepsilon > 0$, one has :
\begin{equation}\label{gooddef}
\bigg| \bigg\{ x \in B \big| |f(x)| < \varepsilon 
 \bigg \} \bigg| \leq
C\left(\displaystyle \frac{\varepsilon}{\sup_{x \in B}|f(x)|}\right)^{\al}|B|.
\end{equation}
\noindent The following elementary properties of $(C, \al)$-good functions will be used.
\begin{enumerate}
\item[(G1)] If $f$ is $(C,\al)$-good on an open set $V$, so is $\lambda
f~\forall~\lambda \in
\bbR$;\\
\item[(G2)] If $f_i, i \in I$ are $(C,\al)$-good on $V$, so is $\sup_{i \in
I}|f_i|$;\\
\item[(G3)] If $f$ is $(C,\al)$-good on $V$ and for some $c_1,c_2\,\textgreater \,0,\, c_1\leq \frac{|f(x)|}{|g(x)|}\leq c_2
\text{ for all }x \in V$, then g is $(C(c_2/c_1)^{\al},\al)$-good on $V$.\\
\item[(G4)] If $f$ is $(C,\al)$-good on $V$, it is $(C',\alpha')$-good
on $V'$ for every $C' \geq C$, $\alpha' \leq \alpha$ and $V'
\subset V$.
\end{enumerate}
One can note that from (G2), it follows that the supremum norm of a vector valued function $\f$ is $(C,\al)$-good 
whenever each of its components is $(C,\al)$-good. Furthermore, in view of (G3), we can replace the norm by an equivalent one, only affecting
$C$ but not $\al$. \\

The next Proposition provides the most important class of good functions.
\begin{proposition}[Lemma 3.2 in \cite{BKM}]\label{goodprop}
Any polynomial $f\in \mathbb{R}[x_1,...,x_s]$ of degree not exceeding $l$ is $(C_{s,l},\frac{1}{sl})$-good on $\mathbb{R}^s$, where 
$C_{s,l}=\frac{2^{s+1}sl(l+1)^{1/l}}{V_s}$. In particular,  constant and linear polynomials are 
$(\frac{2^{s+2}s}{V_s},\frac{1}{s})$-good on $\mathbb{R}^s$.
\end{proposition}
\section{Small Gradients}\label{small-gradients}
For each $t\in \mathbb{Z}_+$, we define $\mathcal{A}_t$ as the set
\[\left \{\bx\in U:\exists p \in \bbZ, \bq \in \bbZ^n \text{ s.t. }
   \left|
\begin{array}{ll}
    |p +  (\bx,  \tilde{\bx}A)\bq| < \displaystyle \frac{\kappa}{2^{nt}}\\
    \|\nabla (\bx,  \tilde{\bx}A) \cdot \bq  \| < \displaystyle \sqrt{\frac{ns}{2r^2}}2^{t/2}\\
    2^t \leq \|\bq\|< 2^{t+1}
\end{array}\right.\right\}\,.\]
It is now immediate that \[\bigcup_{\bq \in \mathbb{Z}^n\setminus \{\mathbf{0}\}}
\mathcal{L}_{small}(\bq)\subseteq \bigcup_{t=0}^{\infty}\mathcal{A}_t\,,\]
since $\forall x\in \mathbb{R}, \psi(x)\leq 1/x$. It is therefore enough to show
\begin{equation}\label{eqn:small estimate}
 \displaystyle \sum_{t=0}^{\infty}|\mathcal{A}_t|<\frac{\xi}{2}\,|U|\,.
\end{equation}
For $\beta\in (0,\frac{1}{2(n+1)})$, we set
\begin{equation}\label{eqn:constants1}
  \delta:=\displaystyle \frac{\kappa}{2^{nt}}, K:=\displaystyle\sqrt{\frac{ns}{2r^2}}\,2^{t/2}, T:=2^{t+1},
\end{equation}  
\begin{equation}
  \varepsilon':=(\delta K T^{n-1})^{\frac{1}{n+1}}=\left(\kappa\,2^{n-1}\,\sqrt{\frac{ns}{2r^2}}\frac{1}{2^{t/2}}\right)^{\frac{1}{n+1}}=
 \left(\frac{\kappa\,2^{n-\frac{3}{2}}\sqrt{ns}}{r}\right)^{\frac{1}{n+1}} \frac{1}{2^{t/2(n+1)}}
 \end{equation}
 and 
 \begin{equation}
 \varepsilon:= 2^{\beta t}\varepsilon'= \left(\frac{\kappa\,2^{n-\frac{3}{2}}\sqrt{ns}}{r}\right)^{\frac{1}{n+1}} \frac{2^{\beta t}}{2^{t/2(n+1)}}.
  \end{equation}

\noindent Define
\begin{equation}\label{unip} u_{\mathbf{x}} : = \left (\begin{array}{rccl}
1 & 0 & \bx& \bx A'+\mathbf{a}_0\\
0 & I_s & I_s& A'\\
0 & 0 & I_{n}\\
\end{array}\right)\end{equation}
\noindent and 
\begin{equation}\label{diag}
g_t := \diag \left(\frac{\varepsilon}{\delta},\frac{\varepsilon}{K},\dots,\frac{\varepsilon}{K},
\frac{\varepsilon}{T},\dots,\frac{\varepsilon}{T}\right)\,.
\end{equation}\\
\noindent Denote by $\La$ the subgroup of $\bbZ^{1+s+n}$
consisting of vectors of the form:
\begin{equation}\label{latt}
\Lambda = \left\{\begin{pmatrix}
p\\
0\\
\vdots\\
0\\
\bq
\end{pmatrix} \mid p \in \bbZ, \bq \in \bbZ^{n} \right\}.
\end{equation}

\noindent It can be easily seen that
\begin{equation}\label{subset}
\mathcal{A}_t \subseteq \tilde{\mathcal{A}}_t:=\{\bx \in U ~:~ \|g_tu_\bx \lambda\| < \varepsilon
~\text{for some}~\lambda \in \La \backslash \{0\}\}.
\end{equation}

\noindent We shall show that if $\kappa$ is taken to be not exceeding $\frac{r}{2^{n-\frac{3}{2}}\sqrt{ns}}$ and $1$ then,
depending on $A$, $\beta$ can be suitably chosen so that 
\begin{equation}\label{smderef1}
|\tilde{\mathcal{A}}_t|\leq \displaystyle K_0\,\kappa^{\frac{1}{s(n+1)}}\,\frac{1}{2^{{\left(\frac{\frac{1}{2(n+1)}-\beta}{s}\right)t}}}\,|U|
\end{equation} for some explicit constant $K_0$ depending on $s,n,U\text{ and }A$ only. One can then set
\begin{equation}\label{eqn:K1} K_1:= \sum_{t=0}^{\infty}\frac{1}{2^{\left(\frac{\frac{1}{2(n+1)}-\beta}{s}\right)t}}\end{equation} and 
reduce $\kappa$ sufficiently to conclude
\[\sum_{t=0}^{\infty}|\mathcal{A}_t|\leq K_0K_1\,\kappa^{\frac{1}{s(n+1)}}|U|<\frac{\xi}{2}\,|U|\,;\] which establishes 
(\ref{eqn:small estimate}).\,\,\,$\Box$\\

The inequality (\ref{subset}) will be proved using the quantitative nondivergence estimate of Kleinbock and Margulis in the next section.
\section{A quantitative nondivergence estimate}
Let $W$ be a finite dimensional real vector space. For a discrete subgroup $\Gamma$ of $W$,  we set $\Gamma_{\bbR}$ to be the 
minimal linear subspace of $W$ containing $\Gamma$. A subgroup $\Ga$ of
$\La$ is said to be primitive in $\La$ if $\Ga = \Ga_{\bbR} \cap
\La$. We denote the set of all nonzero primitive subgroups of $\Gamma$ by $\mathcal{L}(\Gamma)$. Let $j:=\dim (\Gamma_{\bbR})$ be the \emph{rank} of $\Gamma$. We say that 
$\bw\in \bigwedge^j (W)$ represents $\Gamma$ if 
\[\bw=\left \{\begin{array}{rcl}1 &\text{if }j=0\\\bv_1\wedge\cdots\wedge\bv_j &\text{if }j>0\text{ and }\bv_1, \dots,
\bv_j \text{ is a basis of }\Gamma\,.\end{array}\right.\]
In fact, one can easily see that such a representative of $\Gamma$ is always unique up to a sign.\\

A function $\nu:\bigwedge (W)\longrightarrow \bbR_{+}$ is called \emph{submultiplicative} if
\begin{enumerate}
 \item [(i)]$\nu$ is continuous with respect to natural topology on $\bigwedge (W)$;
 \item [(ii)]$\forall t\in \bbR \text{ and }\bw\in \bigwedge (W)$, $\nu(t\bw)=|t|\nu(\bw)$, i.e. it is homogeneous;
 \item [(iii)]$\forall \mathbf{u},\bw \in \bigwedge (W),\nu(\mathbf{u}\wedge\bw)\leq \nu(\mathbf{u})\nu(\bw)$. 
 \end{enumerate}

\noindent In view of property (ii) as given above, without any confusion, we can define $\nu(\Gamma)=\nu(\bw)$, where $\bw$ represents $\Gamma$. \\

Now we shall come to the ``quantitative nondivergence estimate'' which is a generalization of Theorem 5.2 of \cite{KM}.

\begin{theorem}[\cite{BKM}, Theorem $6.2$]\label{BKM1}
Let $W$ be a finite dimensional real vector space, $\Lambda$ a discrete subgroup of $W$ of rank $k$, and a ball
$B = B(x_0,r_0) \subset
\bbR^s$ and a continuous map $H : \tilde{B} \longrightarrow \GL(W)$ be given, where $\tilde{B}=B(x_0,3^kr_0)$.  Take
$C\geq1,\alpha > 0,~0 < \rho < 1$ and $\nu$ be a submultiplicative function $\bigwedge(W)$. Assume that for any $\Gamma \in \mathcal{L}(\Lambda)$,
\begin{enumerate}
\item[(KM1)] the function $x \mapsto \nu(H(x)\Gamma)$ is $(C,\alpha)$-good
on
$\tilde{B}$,
\item[(KM2)]  $\displaystyle \sup_{x\in B}\nu(H(x)\Gamma) \geq
\rho$ and
\item[(KM3)] $\forall~x \in \tilde{B}$, $\# \{ \Gamma \in
\mathcal{L}(\Lambda):\nu( H(x)\Gamma)
 < \rho\} < \infty$.
\end{enumerate}
\noindent Then for every  $\varepsilon'' >0$ one has :
\begin{equation}\label{BKMeq}
|\{x \in B ~:~ \nu(H(x)\mathbf{\lambda}) < \varepsilon'' ~ \text{for
some}~ \mathbf{\lambda} \in \La \backslash \{0\}\}| < k(3^s N_s)^{k}
C \left(\frac{\varepsilon''}{\rho}\right)^{\alpha}|B|.
\end{equation}
\end{theorem}
With the intention of using Theorem \ref{BKM1} to prove (\ref{subset}), we set
$W=\bbR^{1+s+n}$ with basis $\mathbf{e}_0,\mathbf{e}_{*1},\cdots,\mathbf{e}_{*s},\mathbf{e}_1,\cdots,\mathbf{e}_n$, $\Lambda$ as 
given in (\ref{latt}), $B=U$ and $H(\bx)=g_tu_{\bx}$. The submultiplicative function $\nu$ will be chosen, as introduced in \cite[\S7]{BKM}, 
in the following way:\\

Let  $W_*$  be the subspace of $W$ spanned by $\esa, \dots, \esd$. We shall identify $W_{*}^{\perp}$ with $\bbR^{n+1}$ canonically. Also let $\cW$ be the 
ideal of $\bigwedge(W)$ generated by $\bigwedge^{2}(W_*)$, $\pi_{*}$ be the orthogonal 
projection with kernel $\cW$ and $\|\bw\|_e$ be the Euclidean norm of $\pi_{*}(\bw)$. In simple words,
if $\bw$ is written as a sum of exterior products of the base vectors $\mathbf{e}_i$ and $\mathbf{e}_{*i}$,
to compute $\nu(\bw)$, we ignore the components containing exterior products of type $\mathbf{e}_{*i}\wedge\mathbf{e}_{*j},1\leq i\neq j\leq s$, and 
consider the Euclidean norm of rest. It is immediate that $\nu|_W$ agrees with the Euclidean norm. \\

We now seek for proper $C,\al,\rho$ which make (KM1)-(KM3) true. The condition (KM3) can be established for any $\rho\leq1$ exactly in the way it is done in 
\cite[\S7]{BKM}. The following section is devoted to the verification of the remaining ones along with the search for 
the explicit constants.

\section{ Checking (KM1) and (KM2)}\label{ckm1}
\noindent We begin with the explicit computation of $H(\bx)\bw$ for all $\bw\in \bigwedge^k(W_{*}^{\perp})$ and $k=1,\cdots,n+1$. First writing 
$\bx=(x_1,\cdots,x_s)$ and $(\bx,\tilde{\bx}A)=(f_1(\bx),\cdots,f_n(\bx))$, we see that
\begin{enumerate}
 \item $H(\bx)\,\mathbf{e}_0= \frac {\varepsilon}{\delta}\,\mathbf{e}_0$
 \item $H(\bx)\,\mathbf{e}_{*i}= \frac {\varepsilon}{K}\,\mathbf{e}_{*i}\,; \text{ for }1\leq i\leq s$
\item  $H(\bx)\,\mathbf{e}_{i}= \frac {\varepsilon}{\delta}\,f_i(x)\,\eo + \frac {\varepsilon}{K}\,\sum_{j =
1}^{s} \pfij \esj + \frac {\varepsilon}{T}\,\ei\, \text{ for }1\leq i\leq n.$
 \end{enumerate}
 Note that each $f_i(x)$ is a polynomial $x_1,\cdots,x_s$ with degree at most $1$ so that each $\pfij$ is constant. \\
 
 \subsection{Checking (KM1)}: Since $\Lambda=\bbZ^{1+s+n}\cap W_{*}^{\perp}$, any representative $\bw\in \bigwedge^k(W)$ 
 of any subgroup of $\Lambda$ of rank $k$, $1\leq k\leq n+1$, can be written as $\displaystyle \sum_{I}a_{I}\mathbf{e}_{I}$, where each $a_I\in \bbZ$ and 
 $\mathbf{e}_I=\mathbf{e}_{i_i}\wedge \cdots \wedge \mathbf{e}_{i_k}$ with $i_1,\cdots,i_k \in \{0,1,\cdots,n\}, i_1<\cdots<i_k$. \\
 
 Since each component of $\pi_{*}(H(\bx)\bw)$ is a polynomial in $x_1,\cdots,x_s$ 
 with degree at most $1$  in view of 
 (\ref{goodprop}), each of them is $(\frac{2^{s+2}s}{V_s},\frac{1}{s})$-good on $\tilde{U}$. This makes 
 $\|\pi_{*}(H(\bx)\bw)\|$ $(\frac{2^{s+2}s}{V_s},\frac{1}{s})$-good on $\tilde{U}$. As
 \[\displaystyle \frac{1}{2^{\frac{1+s+n}{2}}}\leq \frac{\|\pi_{*}(H(\bx)\bw)\|}{\nu(\pi_{*}(H(\bx)\bw))}\leq 1\,,\]whence, from property 
 (G4) of good functions, $\nu(\pi_{*}(H(\bx)\bw))$ is $(C,\al)$-good with 
 \begin{equation}\label{eqn:C,al}
  C:=\max \left\{\frac{2^{\left(s+2+\frac{1+s+n}{2s}\right)}s}{V_s},1\right\}\text{ and }\al:=\frac{1}{s}\,.
 \end{equation} This verifies (KM1).\,\,\,$\Box$\\

\subsection{Checking (KM2)}\label{ckm2}
Let $\Gamma$ be a subgroup of $\Lambda$ with rank $k$ and $\bw\in \bigwedge^k( W_{*}^{\perp})$ represent $\Gamma$. We first consider the case $k=n+1$. 
So $\bw  = w\, \eo \wedge \ea \wedge \dots \wedge \en$ where $w\in \bbZ\backslash\{0\}$. For any $\bx\in U$, the coefficient of 
$\eo \wedge \esa \wedge \mathbf{e}_2 \wedge \dots \wedge \en$ in $\pi_{*}(H(\bx)\bw)$ is clearly seen to be 
\[w\frac{\varepsilon^{n+1}}{\delta K T^{n-1}}\,.\] Now looking at (\ref{eqn:constants1}), we see that

\begin{equation}\label{eqn:toprank}
 \begin{array}{rcl}\displaystyle \sup_{\bx\in U}\nu(H(\bx)\Gamma)=\displaystyle \sup_{\bx\in U}\nu(H(\bx)\bw)\geq \sup_{\bx\in U}
 \|\pi_{*}(H(\bx)\bw)\| \geq |w\frac{\varepsilon^{n+1}}{\delta K T^{n-1}}|
  \\ \displaystyle =|w|2^{\beta(n+1) t}\frac{(\varepsilon')^{n+1}}{\delta K T^{n-1}}\\ \displaystyle =|w|2^{\beta(n+1) t} \geq \frac{1}{2}.\end{array}
\end{equation}

Assume now $1\leq k\leq n$. To bound the norm of $\|\pi_{*}(H(\bx)\bw)\|$ from below,  we will proceed along the lines of \cite[\S5.3]{G-mani}   
using a technique from \cite{Kleinbock-exponent}. As observed in \cite[\S5.3]{G-mani}, for any $\bx\in U,\,\|\pi_{*}(H(\bx)\bw)\|\geq\|\tilde{g}_{t}\tilde{u}_{\bx}\bw\|$  where
\begin{equation}\label{unipnew} \tilde{u}_{\bx}  = \begin{pmatrix}
1  & \bx&\tilde{\bx}A\\
0  & I_{n}\\
\end{pmatrix},
\end{equation}
\noindent and 
\begin{equation}\label{diag1}
\tilde{g}_{t} =
\diag \left(\frac{\varepsilon}{\delta},
\frac{\varepsilon}{T},\cdots,\frac{\varepsilon}{T}\right).
\end{equation}This inspires us to bound $\sup_{\bx\in U}\|\tilde{g}_{t} \tilde{u}_{\bx} \bw\|$ from below.\\

It follows from ($4.6$) in \cite{Kleinbock-exponent} that
\begin{equation}\label{boundformat1}
\displaystyle \sup_{\bx\in U}\|\tilde{g}_{t} \tilde{u}_{\bx} \bw\|\geq \frac{1}{2^\frac{n+1}{2}}\max
\left\{\left(\frac{\varepsilon^k}{\delta T^{k-1}}\right)\sup_{\bx \in U}\|(\bx,\tilde{\bx}A)
\bc(\bw)\|, \left(\frac{\varepsilon}{T}\right)^k\|\pi(\bw)\|\right\}
\end{equation}
\noindent where $\pi$ is the projection from $\bigwedge(W_{*}^{\perp})$ to $\bigwedge(W_{1\to n})$ and $W_{1\to n}$ stands for the 
span of $\be_1,\cdots,\be_n$ . \\

We recall  that 
\begin{equation}
(\bx,\tilde{\bx}A) = \tilde{\bx}R_A
\end{equation}
\noindent where $R_A$ is defined in (\ref{defR}). Because of this, we can replace in our norm 
calculations, $\sup_{\bx \in U}\|(\bx,\tilde{\bx}A)
\bc(\bw)\|$ by $\sup_{\bx \in U} \|\tilde{\bx}R_A \bc(\bw)\|$. As the functions $1,x_1,\cdots.x_s$ are linearly independent over $\bbR$ on
$U$, the map $\bv \mapsto \displaystyle \sup_{\bx\in U}\|\tilde{\bx}\bv\|$ defines a norm on $(\bigwedge (W_{1\to n}))^{s+1}$ which must be equivalent to the supremum 
norm on $(\bigwedge (W_{1\to n}))^{s+1}$, whence for a constant $K_2>0$ depending on $s,n$ and $U$, we have
\[\sup_{\bx \in U} \|\tilde{\bx}R_A \bc(\bw)\|\geq K_2 \|R_{A}\bc(\bw)\|\,,\] and consequently 
\begin{equation}\label{boundformat2}
\displaystyle \sup_{\bx\in U}\|\tilde{g}_{t} \tilde{u}_{\bx} \bw\|\geq \frac{1}{2^\frac{n+1}{2}}\max
\left\{\left(\frac{\varepsilon^k}{\delta T^{k-1}}\right)K_2 \|R_{A}\bc(\bw)\|, \left(\frac{\varepsilon}{T}\right)^k\|\pi(\bw)\|\right\}\,.
\end{equation}\\
\indent We first note that from Lemma $5.1$ in \cite{Kleinbock-exponent} 
we get that for any $n-s < k \leq n$ and for all but finitely many $\bw \in \bigwedge^{k}(\Lambda)$
\begin{equation}
\|R_A\bc(\bw)\| \geq 1. 
\end{equation}
\noindent It therefore follows that for a constant $K_3>0$ depending alone on $A$, 
\begin{equation}
\displaystyle \sup_{\bx\in U}\|\tilde{g}_{t} \tilde{u}_{\bx} \bw\|\geq \frac{K_2 K_3}{2^\frac{n+1}{2}}\left
(\frac{\varepsilon^k}{\delta T^{k-1}}\right).
\end{equation}
From (\ref{eqn:constants1}) together with the choice $\kappa\leq \frac{r}{2^{n-\frac{3}{2}}\sqrt{ns}}$, 
\begin{equation}\begin{array}{rcl}\displaystyle \frac{\varepsilon^k}{\delta T^{k-1}}= \left(\frac{\kappa\,2^{n-\frac{3}{2}}\sqrt{ns}}{r}\right)^{\frac{k}{n+1}}\times 
\frac{1}{2^{\left(\frac{1}{2(n+1)}-\beta \right)kt}}\times \frac{2^{nt}}{\kappa}\times \frac{1}{2^{(t+1)(k-1)}}\\ \geq \displaystyle
\left(\frac{\kappa\,2^{n-\frac{3}{2}}\sqrt{ns}}{r}\right)\times 
\frac{1}{2^{\left(\frac{1}{2(n+1)}-\beta \right)nt}}\times \frac{2^{nt}}{\kappa}\times \frac{1}{2^{(t+1)(n-1)}}\\ \displaystyle =
\sqrt{\frac{ns}{2r^2}}\times 2^{\left(1-\left(\frac{1}{2(n+1)}-\beta\right)n\right)t}\,.
\end{array}\end{equation}
Picking $\beta \in \left(0,\frac{1}{2(n+1)}\right)$ appropriately, thus we get for all subgroups $\Gamma$ of $\Lambda$ with rank $n-s+1,\cdots,n$, 

\begin{equation}\label{eqn:lb1}
\displaystyle \sup_{\bx\in U}\|\tilde{g}_{t} \tilde{u}_{\bx} \bw\|\geq \frac{K_2 K_3\,\sqrt{ns}}{2^{\frac{n}{2}+1}{r}}
\end{equation} holds true.\\

\indent We will now show how to get analogous lower bounds for subgroups of lower ranks. Recall 
first that a straightforward consequence of (\ref{diocond}) is that we can get constants $\theta, K_4>0$ that depend on $A$ only, with the 
property: for every $1 \leq k \leq n-s$ and $\bw\in \bigwedge^k(\Lambda)$, 

\begin{equation}\label{diocond1}
\|R_{A}
\bc(\bw)\|
\geq K_4 \,\|\pb(\bw)\|^{-\frac {(n - \theta)+1-{k}}{k}}. 
\end{equation}  
\noindent Also for the purposes of obtaining bounds in the lower ranks, we can replace $\pi(\bw)$ in (\ref{boundformat2}) with $\pb(\bw)$ 
as $\|\pi(\bw)\|\geq\|\pb(\bw)\|$. Therefore, for
given $1 \leq k \leq n-s$ and $\bw \in \bigwedge^{k}(\Lambda)$, it suffices to examine
\begin{equation}\label{boundformat3}
\max
\left\{\left(\frac{\varepsilon^k}{\delta T^{k-1}}\right)K_2 K_4 \,\|\pb(\bw)\|^{-\frac {(n - \theta)+1-{k}}{k}}, 
\left(\frac{\varepsilon}{T}\right)^k\|\pb(\bw)\|\right\}.
\end{equation}
\noindent Hence we seek for the solution of equation
\begin{equation}
\frac{K_2K_4}{\delta T^{k-1}}\,y^{-\frac{(n-\theta)+1-k}{k}}=\frac{1}{T^k}\,y
\end{equation}
\noindent which gives $y = (K_2K_4)^{\frac{k}{n-\theta +1}}\left(\frac{T}{\delta}\right)^{\frac{k}{n-\theta+1}}$. This yields (\ref{boundformat3}) is at least
$$ (K_2K_4)^{\frac{k}{n-\theta +1}}\left(\frac{T}{\delta}\right)^{\frac{k}{n-\theta+1}}\left(\frac{\varepsilon}{T}\right)^k$$
$$ =(K_2K_4)^{\frac{k}{n-\theta +1}}\left(\frac{2}{\kappa}\right)^{\frac{k}{n-\theta+1}}2^{\frac{(n+1)kt}{n-\theta+1}}
 \left(\frac{\kappa\,2^{n-\frac{3}{2}}\sqrt{ns}}{r}\right)^{\frac{k}{n+1}} 
\frac{1}{2^{\left(\frac{1}{2(n+1)}-\beta \right)kt}} \frac{1}{2^{(t+1)k}}$$
$$=
(K_2K_4)^{\frac{k}{n-\theta +1}}\left(\frac{2}{\kappa}\right)^{\frac{k}{n-\theta+1}}\frac{1}{2^k}
\left(\frac{\kappa\,2^{n-\frac{3}{2}}\sqrt{ns}}{r}\right)^{\frac{k}{n+1}}2^{\left(\left(\frac{n+1}{n-\theta+1}-1\right)-
\left(\frac{1}{2(n+1)}-\beta \right)\right) kt}$$
$$=
(K_2K_4)^{\frac{k}{n-\theta +1}}\left(\frac{\,2^{n-\frac{3}{2}}\sqrt{ns}}{r}\right)^{\frac{k}{n+1}}
2^{\left(\frac{1}{n-\theta+1}-1\right)k}\kappa^{\frac{-\theta k}{(n-\theta+1)(n+1)}}
2^{\left(\left(\frac{n+1}{n-\theta+1}-1\right)-
\left(\frac{1}{2(n+1)}-\beta \right)\right) kt}.$$
 
Write \[K_5:=\min_{1\leq k\leq n-s}\displaystyle \displaystyle(K_2K_4)^{\frac{k}{n-\theta +1}}\left(\frac{\,2^{n-\frac{3}{2}}\sqrt{ns}}{r}\right)^{\frac{k}{n+1}}
2^{\left(\frac{1}{n-\theta+1}-1\right)k}\,.\] Clearly it depends on $s,n,A \text{ and }U$ only. As $\kappa \leq 1$, by further refining the choice of 
$\beta$ if necessary, one can bound (\ref{boundformat3}) from below by $K_5$; whence, in view of (\ref{boundformat2}), one obtains  
for every $1 \leq k \leq n-s$ and $\bw\in \bigwedge^k(\Lambda)$,
\begin{equation}\label{lb3}
 \displaystyle \sup_{\bx\in U}\|\tilde{g}_{t} \tilde{u}_{\bx} \bw\|\geq \frac{1}{2^\frac{n+1}{2}}\,K_5\,.
\end{equation}
Summarizing the observations (\ref{eqn:toprank}), (\ref{eqn:lb1}) and (\ref{lb3}), we finally confirm (KM2) with the following explicit choice 
of $\rho$
\begin{equation}\label{eqn:rho}
 \displaystyle \min \left \{\frac{1}{2}, \frac{K_2 K_3\,\sqrt{ns}}{2^{\frac{n}{2}+1}{r}},
 \frac{1}{2^\frac{n+1}{2}}\,K_5\right\}\,.
\end{equation}

\section{The proof of (\ref{smderef1})}
The first step towards the proof is the observation 
\[\tilde{\mathcal{A}}_t\subseteq \{x \in U ~:~ \nu(H(x)\mathbf{\lambda}) < \sqrt{1+s+n}\,\varepsilon ~ \text{for
some}~ \mathbf{\lambda} \in \La \backslash \{0\}\}\,.\] This is clear since $\nu|_W$ coincides with the Euclidean norm on $W$. Now applying 
the quantitative nondivergence estimate given by Theorem \ref{BKM1} with $\varepsilon''=\sqrt{1+s+n}\,\varepsilon, C,\al \text{ and }\rho$  as given
in (\ref{eqn:C,al}) and (\ref{eqn:rho}), we have 
\begin{equation}\label{eqn:final}
  |\tilde{\mathcal{A}}_t|\leq |\{x \in U ~:~ \nu(H(x)\mathbf{\lambda}) < \sqrt{1+s+n}\,\varepsilon ~ \text{for
some}~ \mathbf{\lambda} \in \La \backslash \{0\}\}|
\end{equation}
$$\leq (n+1)(3^s N_s)^{n+1}
C (1+s+n)^{\frac{1}{2s}}\left(\frac{\varepsilon}{\rho}\right)^{\frac{1}{s}}|U|$$
$$\leq (n+1)(3^s N_s)^{n+1}
C (1+s+n)^{\frac{1}{2s}}\frac{1}{\rho^{\frac{1}{s}}}\left(\frac{2^{n-\frac{3}{2}}\sqrt{ns}}{r}\right)^{\frac{1}{s(n+1)}}
\kappa^{\frac{1}{s(n+1)}}\frac{1}{2^{\left(\frac{\frac{1}{2(n+1)}-\beta}{s}\right)t}}|U|.$$
Denoting 
$$(n+1)(3^s N_s)^{n+1}
C (1+s+n)^{\frac{1}{2s}}\frac{1}{\rho^{\frac{1}{s}}}\left(\frac{2^{n-\frac{3}{2}}\sqrt{ns}}{r}\right)^{\frac{1}{s(n+1)}}$$
by $K_0$ we hereby conclude (\ref{smderef1}).

\bibliographystyle{amsplain}

\end{document}